\documentclass[12pt,a4paper]{amsart}

\title[Abelian subalgebras from flat tori]{Abelian subalgebras of von Neumann algebras from flat tori in locally symmetric spaces}
\author{Guyan Robertson}
\address{School of Mathematics and Statistics, University of Newcastle, NE1 7RU, U.K.}
\email{a.g.robertson@newcastle.ac.uk}

\usepackage{amsmath}
\usepackage{amscd}

\input{prepictex}
\input{pictex}
\input{postpictex}
\chardef\bslash=`\\ 

\makeatletter
\def\verbatim{\interlinepenalty\@M \@verbatim
  \leftskip\@totalleftmargin\advance\leftskip2pc
  \frenchspacing\@vobeyspaces \@xverbatim}
\makeatother
\newtheorem{theorem}{Theorem}[section]
\newtheorem{corollary}[theorem]{Corollary}
\newtheorem{lemma}[theorem]{Lemma}
\newtheorem{proposition}[theorem]{Proposition}

\theoremstyle{definition}

\newtheorem{remark}[theorem]{Remark}

\newcommand{\VN}{{\text{\rm VN}}}

\newcommand{\vp}{\varepsilon}
\newcommand{\cl}[1]{{\mathcal{#1}}}
\newcommand{\bb}[1]{{\mathbb{#1}}}

\newcounter{picture}

\DeclareMathOperator{\supp}{supp}

\newcommand{\G}{{\Gamma}}

\newcommand{\SL}{{\text{\rm{SL}}}}

\newcommand{\tto}{{\text{\rm{II}}}_{1}}

\newcommand{\SO}{{\text{\rm{SO}}}}

\newcommand{\vn}{{\text{\rm VN}}}
\newcommand{\hd}{{\text{\rm hd}}}

\begin{document}
\def\Proof  {{\bf Proof.}\par}
\def\bs {\backslash}
\def\diam   {{\text{\rm diam}}}

\begin{abstract}
Consider a compact locally symmetric space $M$ of rank $r$, with fundamental group $\Gamma$.
The von Neumann algebra $\vn(\Gamma)$ is the convolution algebra of functions $f\in\ell_2(\Gamma)$ which act by left convolution on $\ell_2(\Gamma)$.  Let $T^r$ be a totally geodesic flat torus of dimension $r$ in $M$ and let $\Gamma_0\cong\bb Z^r$ be the image of the fundamental group of $T^r$ in $\Gamma$. Then $\vn(\Gamma_0)$ is a maximal abelian $\star$-subalgebra of $\vn(\Gamma)$ and its unitary normalizer is as small as possible.
If $M$ has constant negative curvature then the Puk\'anszky invariant of $\vn(\Gamma_0)$ is $\infty$.
\end{abstract}

\maketitle

\section{Introduction}

If $\Gamma$ is a group, then the von Neumann algebra $\vn(\Gamma)$ is the convolution algebra
$$\vn(\Gamma)=\{f\in\ell^2(\Gamma)\, :\, f\star \ell^2(\Gamma)\subseteq\ell^2(\Gamma)\}\,.$$
It is well known that if $\Gamma$ is an infinite conjugacy class [ICC] group then $\vn(\Gamma)$ is a factor of type $\tto$.
If $\Gamma_0$ is a subgroup of $\Gamma$, then
$\vn(\Gamma_0)$ embeds as a subalgebra of $\vn(\Gamma)$ via $f\mapsto \overline f$, where
\begin{equation*} \overline f(x) =\begin{cases}
f(x)& \text{if $x\in\Gamma_0$},\\
0& \text{otherwise}.
\end{cases}
\end{equation*}
This article is concerned with examples where $\Gamma_0$ is an abelian subgroup
of $\Gamma$ and $\vn(\Gamma_0)$ is a maximal abelian $\star$-subalgebra (masa) of $\vn(\Gamma)$.

If $\cl A$ is a von Neumann subalgebra of a von Neumann
algebra $\cl M$ then the {\it unitary normalizer} $N(\cl A)$  is the set of unitaries $u$
in $\cl M$ such that $u\cl Au^{-1} = \cl A$.
The subalgebra $\cl A$ is said to be {\it singular} if $N(\cl A)
\subseteq \cl A$, so that the only normalizing unitaries already
belong to $\cl A$.

Let $\Gamma$ be a torsion free cocompact lattice in a semisimple Lie group $G$ of rank $r$ with no centre and no compact factors.  Consider the Riemannian symmetric space  $X=G/K$ and the compact locally symmetric space $M=\Gamma\backslash X$.
A {\it flat} in $X$ is an isometrically embedded euclidean space in $X$.
The {\it rank} $r$ of $X$ is the dimension of a maximal flat in $X$.

Suppose that $T^r$ is a totally geodesic flat torus of dimension $r$ in $M$. Let $\Gamma_0\cong\bb Z^r$ be the image of the fundamental group $\pi(T^r)$ under the natural monomorphism from  $\pi(T^r)$ into $\Gamma=\pi(M)$. We show that $\vn(\Gamma_0)$ is a singular masa of $\vn(\Gamma)$. In fact we prove two stronger results :  Theorems \ref{a4o} and \ref{a7} below. Fix $g_1\in \Gamma_0$ and let $\Gamma_1\cong\bb Z$ be the subgroup of $\Gamma_0$ generated by $g_1$. Then we have inclusions
$$\vn(\Gamma_1)\subseteq \vn(\Gamma_0)\subseteq \vn(\Gamma).$$
The result below implies that if $g_1$ is the homotopy class of a regular geodesic (as defined subsequently) then $\vn(\Gamma_0)$ is the {\it unique} masa of $\vn(\Gamma)$
containing $\vn(\Gamma_1)$.

\begin{theorem}\label{a4o}
Let $g_1\in\Gamma_0$ be the class of a regular closed geodesic $c$ in $T^r$, and let $\Gamma_1\cong\bb Z$ be the subgroup of $\Gamma_0$ generated by $g_1$.
Let $u$ be a unitary operator in $\cl M$ such that $u \vn(\Gamma_1) u^{-1} \subseteq \vn(\Gamma_0)$.
Then $u\in \vn(\Gamma_0)$.
\end{theorem}

The second main result (Theorem \ref{a7}) implies that $\vn(\Gamma_0)$ is a {\it strongly singular} masa in the sense of  \cite{SS}. This improves a result of \cite[Theorem 4.9]{RSS}, which proved strong singularity under the additional hypothesis that the diameter of  $T^r$ is small. The two new ideas leading to this improvement are the use of the Amenable Subgroup Theorem (Lemma \ref{SST})
and  the replacement of the Furstenberg Boundary by the Tits Boundary in the subsequent arguments.

 If $M$ has constant negative curvature then it is proved in Theorem \ref{puk} that the Puk\'anszky invariant of $\vn(\Gamma_0)$ is $\infty$.

\section{Preliminaries}

We first recall some concepts which are needed for the statements and proofs of the results. Let $G$ be a semisimple Lie group with no centre
and no compact factors.
Let $X=G/K$ be the associated symmetric space, where $K$ is a maximal compact
subgroup of $G$. Then $X$ is a contractible space of nonpositive curvature.
A geodesic $L$ in $X$ is called {\it regular} if it lies in only one
maximal flat; it is called {\it singular} if it is not regular.
Let $F$ be a maximal flat in $X$ and let $p\in F$.  Let $S_p$ denote the
union of all the singular geodesics through $p$.  Then $F-S_p$ has finitely many connected
components,  called {\it Weyl chambers} with origin $p$.

If $A$, $B$ are subsets of $X$, and $\delta>0$, then the notation
$A\displaystyle\operatornamewithlimits{\subset}_{\delta}B$ means that $d(a,B)\leq\delta$,
for all $a\in A$.
Define the Hausdorff distance between $A$ and $B$ to be
\begin{equation}\label{eq4.5}
\hd(A,B)=\inf\{\delta\leq\infty\, :\, A\underset{\delta}{\subset}B\text{ and
}B\underset{\delta}{\subset} A\}\,.
\end{equation}
Any complete geodesic $L$ in $X$ is the union of two geodesic rays which intersect
at their common origin.
Define an equivalence relation $\sim$ on the set of geodesic rays in $X$ by
\begin{equation}\label{equiv}
L_1\sim L_2\Longleftrightarrow\hd(L_1,L_2)<\infty\,.
\end{equation}

The sphere at infinity $X(\infty)$ is the set of
equivalence classes of geodesic rays in $X$ \cite[Chapter III]{GJT}.
Denote by $L(\infty)$ the class in $X(\infty)$ of a geodesic ray $L$.
The set $X(\infty)$ may be given the structure of a spherical building whose maximal simplices are called Weyl chambers at infinity, and there is a natural action of $G$ on $X(\infty)$.
The parabolic subgroups of $G$ are the stabilizers $G_z=\{x\in G \, :\, x(z)=z\}$, for some
$z\in X(\infty)$ \cite[Proposition 3.8]{GJT}. Moreover, $G_z$ is a minimal parabolic subgroup of $G$
if and only if $z=L(\infty)$, where $L$ is a geodesic ray in a Weyl chamber
\cite[pp 248--9]{BGS}.
Two such minimal parabolic subgroups $G_{z_1}$, $G_{z_2}$ coincide if and only if $z_1, z_2$
belong to the same Weyl chamber at infinity \cite[Proposition 3.16]{GJT}.

If $F$ is a maximal flat in $X$, then the restriction of the equivalence
relation $\sim$ to rays in $F$ allows one to define the sphere at infinity $F(\infty)$. There is a natural embedding of $F(\infty)$ into $X(\infty)$, and it is convenient to identify $F(\infty)$ with its image in $X(\infty)$.

Now let $\Gamma$ be a torsion free cocompact lattice in $G$.  Then $\Gamma$  acts
freely
on the symmetric space $X=G/K$ and the quotient manifold $M=\Gamma\backslash
X$
has universal covering space $X$.  The manifold $M$ is a compact locally symmetric
space of nonpositive curvature, with fundamental group $\Gamma$,
and $\Gamma$ acts freely on $X$.

Let $T^r\subset M$ be a totally geodesic embedding of a flat $r$-torus
in $M$. Choose and fix a point $\xi\in T^r$. Consider the fundamental groups $\Gamma=\pi(M,\xi)$ and $\Gamma_0=\pi(T^r,\xi)\cong\bb Z^r$.
Since no geodesic loop in $M$ is null-homotopic,  the inclusion $i:T^r\to M$ induces an injective homomorphism $i_*:\Gamma_0\to\Gamma$.
We identify $\Gamma_0$ with its image in $\Gamma$.
There is a  maximal flat $F_0\cong \bb R^r$ in $X$ such that $\Gamma_0$ acts cocompactly by translations upon $F_0$, and $p(F_0)=T^r$
\cite[Theorem II.7.1]{BH}.

The flat $F_0$ is the unique $\Gamma_0$-invariant flat in $X$. For if $F_1$ is another such, then since $\Gamma_0$ acts isometrically and the action on
$F_0$ is cocompact, we have $F_0\underset{\delta}{\subset}F_1$, for some $\delta>0$. Therefore $F_0=F_1$, by \cite[Lemma 7.3 (iv)]{Mos}, applied to a maximal flat containing $F_1$.

Choose $\tilde\xi\in F_0$ such that $p(\tilde\xi)=\xi$.  For any element $x\in\Gamma$, there is a unique geodesic loop based at $\xi$ which represents $x$.  This is the loop $c$ of shortest length in the class $x$ and it is the projection of the geodesic segment $[\tilde\xi,x\tilde\xi]$ in $X$.

The following technical lemma will play a crucial role later, in our improvement to \cite[Theorem 4.9]{RSS}.
\begin{lemma}\label{SST}
Let $z\in F_0(\infty)$ lie in a Weyl chamber at infinity.
Then
$$G_z\cap\Gamma=\Gamma_0.$$
\end{lemma}

\begin{proof} The group $G_z=\{x\in G \, :\, \, x(z)=z\}$ is a minimal parabolic subgroup of $G$,
and so has a cocompact solvable normal subgroup.
Therefore $G_z$ is amenable, as is the discrete subgroup $\Gamma_z= G_z\cap\Gamma$.
Also $\Gamma_z \supseteq \Gamma_0$, since $\Gamma_0$ acts by translation upon $F_0$
and so stabilizes each point of $F_0(\infty)$.

By the Amenable Subgroup Theorem \cite[Corollary B]{AB}, there is a $\Gamma_z$-invariant flat $F_z$ in $X$.  Since $\Gamma_z \supseteq \Gamma_0$, $F_z$ is $\Gamma_0$-invariant and so $F_z=F_0$, by the remarks preceding this lemma.
Thus $\Gamma_zF_0=F_0$.

If $x\in\Gamma_z$ then $x\tilde\xi\in F_0$, since $\tilde\xi\in F_0$.  Therefore the geodesic segment $[\tilde\xi,x\tilde\xi]$ in $F_0$ projects to a closed geodesic in $T^r$, whose class in $\Gamma$ is precisely $x$. Therefore $x\in\Gamma_0$.
\end{proof}

\section{Singularity results for abelian subalgebras}\label{3}

The terminology and notation introduced in the previous section will be used without further comment.
A {\it regular} geodesic in $M$ is, by definition, the image of a regular geodesic under the covering projection $p:X\to M$. It follows from \cite[\S 11]{Mos} that $T^r$ contains a closed regular geodesic.
We now prove Theorem \ref{a4o}, which we restate here, for convenience.

\begin{theorem}\label{a4}
Let $g_1\in\Gamma_0$ be the class of a regular closed geodesic $c$ in $T^r$, and let $\Gamma_1\cong\bb Z$ be the subgroup of $\Gamma_0$ generated by $g_1$.
Let $u$ be a unitary operator in $\cl M$ such that $u \vn(\Gamma_1) u^{-1} \subseteq \vn(\Gamma_0)$.
Then $u\in \vn(\Gamma_0)$.
\end{theorem}

\begin{proof}
Suppose that $x_0 \in \supp\, u$. We must prove that $x_0 \in \G_0$.

Since $u\in \ell^2(\Gamma)$, there are only a finite number of cosets $\G_0y$, with $y\in \Gamma$, such that
$\|u|_{\G_0y}\|_2 \ge |u(x_0)|$.  Call these cosets $\G _0y_1, \dots ,\G
_0 y_n$.  We claim that
\begin{equation}\label{modified}
 x_0\G_1 \subset\G _0 y_1 \cup \dots \cup \G _0y_n.
\end{equation}  To prove this, note that if $z \in \G_1$ then $u
\star \delta_z \star u^{-1} = f$ is a unitary operator which lies in $\vn(\G_0)$, by hypothesis. Therefore
\begin{align*}
|u(x_0)|
&= |(u \star \delta_z)(x_0z)| \\
&= |(f \star u)(x_0 z)| \\
&=\left| \sum_{t\in \Gamma_0} f(t)u(t^{-1}x_0z)\right| & \text{(since $\supp\, f \subseteq \Gamma_0$)} \\
&\le \|u|_{\G_0x_0z}\|_2  & \text{(since $f$ is unitary)}.
\notag
\end{align*}
This shows that $x_0z \in \G _0y_1 \cup \dots \cup \G _0y_n$,
which proves (\ref{modified}).

We now show that (\ref{modified}) implies that $x_0 \in \G_0$.
Lift $c$ to a regular geodesic $L$ in $X$ through a point $\tilde\xi\in X$.
Regularity means that $L$ lies in a {\it unique} maximal flat $F_0$
and $p(F_0)=T^r$. The diagram below illustrates the case
$X=\SL_3(\bb R)/\SO_3(\bb R)$, where $r=2$ and there are six Weyl chambers in
the flat $F_0$ with  origin $\tilde\xi$.

\begin{figure}[htbp]
\centerline{
\beginpicture
\setcoordinatesystem units <1cm, 1.732cm>
\setplotarea  x from -2 to 2,  y from -4 to 1.5
{\tiny
\putrule from  0 -1.3 to 0 -2
\arrow <10pt> [.2, .67] from  0 -1.8 to  0 -2
\plot  1.5 -0.8  0 0  -1.5 0.8 /
\put{$L$} [l] at 1.8 -0.8
\put {$_{\bullet}$} at 0 0
\put {$\tilde\xi$} [b,l] at 0.27 0.05
\put{
\beginpicture
\setcoordinatesystem units  <0.5in, 0.25in>   
\setplotarea x from 0 to 2, y from -2 to 2  
\circulararc 140 degrees from 0.7 0.2 center at  1 0.4         
\circulararc 72 degrees from 1.2 -0.05 center at  1 -0.6
\ellipticalarc axes ratio  2:1  360 degrees  from  2.5 0
center at 1.1  0
\setquadratic  
\plot
0.9 -0.25    0.8  -0.8    0.9  -1.4  /
\put {$_c$}  at  0.7 -0.8
\put {$T^r$}  at  2.9 0
\endpicture
}
at 0.6 -3
\setdashes
\putrule from -2 0 to 2 0
\plot  -1 -1  0 0  1 1  /
\plot  1 -1  0 0  -1 1 /
}  
\endpicture
}
\end{figure}

Let $P=[\tilde\xi,g_1\tilde\xi]\subseteq F_0$, so that
$\Gamma_1 P=L$, since $g_1$ acts on $L$ by translation.
Let $\delta=\max\{d(y_jp,p)\, :\,1\leq j\leq n,\,p\in P\}$.
Then for $1\le j\le n$, we have $y_j P\underset{\delta}{\subset} P$, and so
\begin{equation}
\Gamma_0y_j P\underset{\delta}{\subset}\Gamma_0 P\subset F_0\,.
\end{equation}

It follows from (\ref{modified}) that
\begin{equation}
x_0L=x_0\Gamma_1P\underset{\delta}{\subset}F_0\,.
\end{equation}
Since $x_0L$ is a regular geodesic, \cite[Lemma 7.3(iii)]{Mos} implies that $d(x_0L, F_0)=0$.
Consequently $x_0L\subset F_0$, by \cite[Lemma 3.7]{Mos}.
In particular, $x_0\tilde\xi\in F_0$.
Therefore the geodesic segment $[\tilde\xi,x_0\tilde\xi]$ projects to a closed geodesic in $T^r$ whose class in $\Gamma$ is precisely $x_0$.  Hence $x_0\in\Gamma_0$.
\end{proof}

\begin{remark}
Theorem \ref{a4} implies that $\vn(\Gamma_0)$ is a singular masa of $\vn(\Gamma)$ and that it is also the unique
masa of $\vn(\Gamma)$ containing $\vn(\Gamma_1)$.
Since closed geodesics are dense in the set of all geodesics of $T^r$,
we can choose regular geodesics $c_1, c_2, \dots , c_r$ in $T^r$ which lift to regular geodesics
in linearly independent directions in $F_0$.
Applying Theorem \ref{a4} to each of the geodesics $c_j$ shows that
$\vn(\Gamma)$ contains a masa $\cl A =\vn(\Gamma_0)$ with the following property :

$\cl A$ contains abelian subalgebras $\cl B_j$, $1\le j\le r$ such that

\begin{itemize}
\item[(a)] $\cl A$ is the {\it unique} masa containing $\cl B_j$\,;
\item[(b)] $\cl A$ is generated by $\cl B_1\cup\cl B_2\cup\dots\cup \cl B_r$\,;
\item[(c)] $\cl B_i$, $\cl B_j$ are {\it orthogonal} for $i\ne j$, in the sense that
$Tr(b_ib_j)=0$ whenever  $Tr(b_i)=Tr(b_j)=0$ \cite[Definition 2.2]{Po}.
\end{itemize}
This construction shows how one can see something of the rank of
the locally symmetric space $M$ in the group von Neumann algebra of the fundamental group  $\Gamma=\pi(M)$. This is of interest in connection with a conjecture of A. Connes.

\medskip

\noindent {\bf Rigidity Conjecture}.  If ICC groups $\Gamma_1$,$\Gamma_2$ have Property (T)
of Kazhdan, then 
\[
\VN(\Gamma_1)\cong\VN(\Gamma_2) \Rightarrow \Gamma_1\cong \Gamma_2.
\]

In our setup, one consequence of the truth of Connes' conjecture would be that the rank of $M$ is determined by $\VN(\Gamma)$.
\end{remark}

Now we define a relative version of the notion of a strongly singular masa
defined in \cite{SS}.
Let $\cl A \subseteq \cl C \subseteq \cl M$ be von Neumann subalgebras of a type $\tto$
factor $\cl M$, and let ${\mathbb{E}}_{\cl N}$ denote the unique trace preserving
conditional expectation onto any von Neumann subalgebra $\cl N$ of
$\cl M$. Say that $\cl A \subseteq \cl C$ is  {\it {a strongly singular pair}} of von Neumann subalgebras of $\cl M$ if for all von Neumann subalgebras $\cl B$ with $\cl A \subseteq \cl B \subseteq \cl C$ the inequality
\begin{equation}\label{eqa}
\|\bb E_{\cl B} - \bb E_{u\cl Bu^*}\|_{\infty,2} \ge \|(I-\bb E_{\cl C})(u)\|_2
\end{equation}
holds for all unitaries $u\in\cl M$. [As in \cite{RSS}, the notation $\|\ T \|_{\infty,2}$ means that the norm of the linear map $T$ is taken relative to operator norm on its domain and the $\ell^2$ norm on its range.
If $\cl A=\cl C$ then this reduces to the definition of a strongly singular subalgebra given in \cite{SS}.

The next two results are mild generalizations of \cite[Lemma 2.1]{RSS} and \cite[Lemma 4.1]{RSS}. The proofs are included for completeness.

\begin{lemma}\label{a1}
Let $\cl A \subseteq \cl C \subseteq \cl M$ be von Neumann subalgebras of a type $\tto$ factor $\cl M$.
Suppose that, given $\vp>0$ and a unitary operator $u\in\cl M$, there exists a unitary $v\in \cl A$, such that
\begin{equation}\label{eqb}
\|\bb E_{\cl C}(u^*v u) - \bb E_{\cl C}(u^*)v \bb E_{\cl C}(u)\|_2 < \vp.
\end{equation}
Then $\cl A \subseteq \cl C$ is a strongly singular pair.
That is, (\ref{eqa}) holds, whenever $\cl A \subseteq \cl B \subseteq \cl C$.
\end{lemma}

\begin{proof}
We have
\begin{align}
\|\bb E_{\cl B} - \bb E_{u\cl Bu^*}\|^2_{\infty,2} &\ge \|v-\bb E_{u\cl
Bu^*}(v)\|^2_2 \qquad\qquad\qquad \text{since $v\in\cl B$}\\
&= \|v-u\bb E_{\cl B}(u^*vu)u^*\|^2_2\nonumber\\
&= \|u^*vu-\bb E_{\cl B}(u^*vu)\|^2_2\nonumber\\
&= 1 - \|\bb E_{\cl B}(u^*vu)\|^2_2 \qquad\qquad \text{by orthogonality}\nonumber\\
&\ge 1 - \|\bb E_{\cl C}(u^*vu)\|^2_2 \qquad\qquad\qquad\qquad \nonumber\\
&\ge 1 - (\|\bb E_{\cl C}(u^*)v\bb E_{\cl C}(u)\|_2 + \vp)^2\qquad\quad\, \text{by (\ref{eqb})}\nonumber\\
&\ge 1-(\|\bb E_{\cl C}(u)\|_2 + \vp)^2 \nonumber\\
\label{eq2.6}
&= \|(I-\bb E_{\cl C})(u)\|^2_2 - \vp^2-2\vp\|\bb E_{\cl C}(u)\|_2\,.
\end{align}
Since $\vp>0$ was arbitrary, the result follows.
\end{proof}

\begin{lemma}\label{a2}
Let $\Gamma_1 < \Gamma_0 < \Gamma$ be subgroups of a discrete ${\mathrm{I.C.C.}}$ group.
The following condition implies that
$\vn(\Gamma_1) \subseteq \vn(\Gamma_0)$ is  {\it {a strongly singular pair}}. \\
If $x_1,\ldots,x_m \in \Gamma$ and
\begin{equation}\label{eq6}
\Gamma_1\subseteq\bigcup_{i,j}x_i\Gamma_0 x_j\, ,
\end{equation}
then $x_i\in \Gamma_0$ for some $i$.
\end{lemma}

\begin{proof}
The condition in question is equivalent to the following: \\
If
$x_1,\ldots ,x_n,y_1,\ldots ,y_n \in \Gamma \backslash \Gamma_0$, then there
exists $\gamma_0\in \Gamma_1$ such that
\begin{equation}\label{eqf}
x_i\gamma_0y_j\notin \Gamma_0,\ \ 1 \leq i, j \leq n.
\end{equation}
To see this replace $x_i$ by $x_i^{-1}$, replace $y_j$ by $y_j^{-1}$
and replace each of the sets $\{x_1,\ldots, x_n\}$, $\{y_1,\ldots,
y_n\}$ by their union, which is renamed $\{x_1,\ldots, x_m\}$.

Let $\cl A= \vn(\Gamma_1)$ and $\cl C = \vn(\Gamma_0)$.
Given $\vp>0$ and a unitary operator $u\in\vn(\Gamma)$, we must show that there exists a unitary $v\in \cl A$, such that (\ref{eqb}) is satisfied.
To do this, approximate $u$ by a finite linear combination of group elements $y_1, \dots y_n, y_{n+1}, \dots , y_p$, where $y_1, \dots y_n \not\in\Gamma_0$ and $y_{n+1}, \dots , y_p \in \Gamma_0$. Let $x_i=y_i^{-1}, 1\le i\le p$ and choose $\gamma_0\in\Gamma_1$ satisfying
(\ref{eqf}). Now
\begin{equation*}
\bb E_{\cl C}(x_i\gamma_0 y_j) = \bb E_{\cl C}(x_i)\gamma_0 \bb E_{\cl C}(y_j),
\qquad 1\le i, j \le n,
\end{equation*}
both sides being zero if $i\le n$ or $j\le n$, since $x_i\gamma_0 y_j$ is then orthogonal to $\cl B$. The equation (\ref{eqb}) follows by taking a close enough approximation.
\end{proof}

The next result implies that $\vn(\Gamma_0)$ is a strongly singular masa of $\vn(\Gamma)$ in the sense of \cite{SS}. It improves \cite[Theorem 4.9]{RSS}, by removing a superfluous hypothesis on the diameter of the embedded torus $T^r$.

\begin{theorem}\label{a7}
Let $g_1$ be the class of a regular closed geodesic $c$ in $T^r$, and let $\Gamma_1\cong\bb Z$ be the subgroup of $\Gamma_0=\pi(T^r)$ generated by $g_1$.
Then $\vn(\Gamma_1) \subseteq \vn(\Gamma_0)$ is a strongly singular pair.

\end{theorem}

\begin{remark} Theorem \ref{a7} is clearly closely related to Theorem \ref{a4}, but neither result
appears to contain the other.
\end{remark}

\begin{proof}
Lift $c$ to a regular geodesic $L$ in $X$ through $\tilde\xi$, where  $p(\tilde\xi)=\xi$
and $L=L^+\cup L^-$ is a union of two geodesic rays with common origin $\tilde\xi$.
Regularity means that $L$ lies in a {\it unique} maximal flat $F_0$
and $p(F_0)=T^r$.

Suppose that (\ref{eq6}) holds. That is, we have elements $x_1,\ldots,x_m \in\Gamma$ such that
\begin{equation}\label{eq6.1}
\Gamma_1\subseteq\bigcup_{i,j}x_i\Gamma_0 x_j\,.
\end{equation}
Now $g_1$ acts on $L$ by translation. Let $P=[\tilde\xi,g_1\tilde\xi]$, so that $L=\Gamma_1P$.
Let $\delta=\max\{d(x_jp,p)\, :\,1\leq j\leq n,\,p\in P\}$.
For $1\leq j\leq n$, this implies that $x_j P\underset{\delta}{\subset} P$ and
so
\begin{equation}
\Gamma_0x_j P\underset{\delta}{\subset}\Gamma_0 P\subset F_0\,.
\end{equation}
Hence, for each $i$, $j$, we have
$x_i\Gamma_0x_jP\underset{\delta}{\subset}x_iF_0$.
It follows from (\ref{eq6.1}) that
\begin{equation}
L=\Gamma_1P\underset{\delta}{\subset}x_1F_0\cup x_2F_0\cup\cdots\cup x_mF_0\,.
\end{equation}
Let $z=L(\infty)\in F_0(\infty)$.
Now for $1\le j\le m$, the element $x_jg_1x_j^{-1}$ acts by translation upon the maximal flat $x_jF_0$ and hence preserves each boundary point of $x_jF_0$.
Express each $x_jF_0$ as a (finite) union of Weyl chambers $W_{\alpha,j}$ with base vertex
$x_j\tilde\xi$. Thus
\begin{equation}\label{subsec}
L^+\underset{\delta}{\subset}\bigcup_{\alpha, j} W_{\alpha,j}\,.
\end{equation}
Now for each $\alpha, j$ the function
$p\mapsto d(p, W_{\alpha,j})$ is convex on $L$ by \cite[Lemma 3.6]{Mos}.
According to (\ref{subsec}), we have
\begin{equation}
\min_{\alpha, j} d(p, W_{\alpha,j})\le \delta
\end{equation}
for all $p\in L$.
This implies that for some $\alpha, j$, the function $p\mapsto d(p, W_{\alpha,j})$
is monotonically decreasing on $L^+$.  Choose such $\alpha, j$. Then
$$L^+\underset{\epsilon}{\subset}W_{\alpha,j}$$ for some $\epsilon >0$.
By \cite[Lemma 7.3(i)]{Mos}, there is a geodesic ray $L'\in x_jF_0$ which is asymptotic to $L$.
Thus $z=L^+(\infty)=L'(\infty)\in x_jF_0(\infty)$ and $x_jg_1x_j^{-1}z=z$.

Since $L$ is a regular geodesic, $z$ lies in a Weyl chamber at infinity and it follows from Lemma \ref{SST} that  $x_jg_1x_j^{-1}=h_0\in \Gamma_0$.
Since $x_jg_1^nx_j^{-1}=h_0^n$, we have $x_jg_1^n=h_0^nx_j$. Therefore
$$d(x_jg_1^n\tilde\xi, h_0^n\tilde\xi)=
d(h_0^nx_j\tilde\xi, h_0^n\tilde\xi)=
d(x_j\tilde\xi,\tilde\xi)\,.$$
Thus $d(x_jL, F_0)<\infty$, from which it follows that  $x_jL\subset F_0$, by \cite[Lemma 7.3(iii) and Lemma 3.7]{Mos},
applied to $x_jL^\pm$.
In particular, $x_j\tilde\xi\in F_0$.
Hence the geodesic segment $[\tilde\xi,x_j\tilde\xi]$ projects to a closed geodesic in $T^r$ whose class in $\Gamma$ is precisely $x_j$.  It follows that $x_j\in\Gamma_0$.
\end{proof}

\bigskip

\section{The Puk\'anszky invariant in constant negative curvature}

Let $\Gamma_0$ be an abelian subgroup of a countable group $\Gamma$ such that $\cl  A =\vn(\Gamma_0)$ is a masa
of $\cl  M = \vn(\Gamma)$. Recall that $\cl  A$ is the von Neumann subalgebra of $B(\ell^2(\Gamma))$ defined by the left convolution operators
$$\lambda(f) : \phi\mapsto f\star \phi$$
where $f\in\ell^2(\Gamma_0)$ and $f\star \ell^2(\Gamma)\subseteq\ell^2(\Gamma)$.
The algebra $\cl A$ also acts on $\ell^2(\Gamma)$ by right convolution 
$$\rho(f) : \phi\mapsto \phi\star f ,$$
where $f\in \cl  A$. Let $\cl  A^{opp}$ be the von Neumann subalgebra of $B(\ell^2(\Gamma))$
defined by this right action of $\cl A$.

Let ${\mathcal B}$ be the von Neumann algebra generated by $\cl  A \cup \cl  A^{opp}$
and let $p$ denote the orthogonal projection of $\ell^2(\Gamma)$ onto the closed
subspace generated by $\cl  A$. Then $p$ is in the centre
of ${\mathcal B}'$ and ${\mathcal B}'p$ is abelian.
The von Neumann algebra ${\mathcal B}'(1 - p)$ is of type $I$ and may therefore be expressed as a direct sum $\mathcal B_{n_1}\oplus \mathcal B_{n_2} \oplus \dots$ of algebras
$\mathcal B_{n_i}$ of type $I_{n_i}$, where $1\le n_1 < n_2 < \dots \le \infty$. The {\it Puk\'anszky invariant} \cite {po2} is the set $\{n_1, n_2, \dots \}$.  It is an isomorphism invariant of the
pair $(\cl  A,\cl  M)$, since any automorphism of $\cl  M$ is implemented by a unitary
in $B(\ell^2(\Gamma))$. It has been shown \cite[Corollary 3.3]{NS} that all subsets of the natural numbers can be realized as the Puk\'anszky invariant of some masa of the hyperfinite $II_1$ factor.

A subgroup $\Gamma_0$ of a group $\Gamma$ is {\it malnormal}
if $g^{-1}\Gamma_0g \cap \Gamma_0 = \{ 1 \}$ for all
$g \in \Gamma - \Gamma_0$.
Recall the following result from \cite[Proposition 3.6]{RS},
which we shall apply in proving Theorem \ref{puk}.
\begin{proposition} \label{3e} Suppose that $\Gamma_0$ is an abelian subgroup of a countable group $\Gamma$
such that  $\cl  A = \vn(\Gamma_0)$ is a masa of $\vn(\Gamma)$.
If $\Gamma_0$ is malnormal then the Puk\'anszky invariant
of $\cl A$ is $n = \#(\Gamma_0 \backslash \Gamma /\Gamma_0)$.
\end{proposition}
Return now to the setup of Theorem \ref{a4}. Thus $T^r$ is a totally geodesic flat torus of dimension $r$ in the compact locally symmetric space $M$ and $\Gamma_0\cong\bb Z^r$ is the image of the fundamental group $\pi(T^r,\xi)$ in $\Gamma$.
The maximal flat $F_0$ in $X$ covers $T^r$ and the element $\tilde\xi\in F_0$ projects to $\xi\in M$.

It is not always true that $\Gamma_0$ is malnormal in $\Gamma$.
Nevertheless, there is a weaker result.
\begin{lemma}\label{conjr}
Suppose that $g\in\Gamma$ and that $g^{-1}\Gamma_0g \cap \Gamma_0$ contains an element
$x_0\not= 1$ which is the class of a regular closed geodesic $c$ in $T^r$. Then $g\in \Gamma_0$.
\end{lemma}

\begin{proof}
This follows immediately from Theorem \ref{a4}.
\end{proof}

\begin{corollary}\label{conjrcor}
Suppose that $g\in\Gamma$ and that $g^{-1}\Gamma_0g \cap \Gamma_0$ contains a free abelian group of rank $r$. Then $g\in \Gamma_0$.
\end{corollary}

\begin{proof}
Combine Lemma \ref{conjr} and  \cite[Lemma 11.1]{Mos}.
\end{proof}

\begin{corollary}\label{mal}
Suppose that $M$ has strictly negative curvature  and that $x_0$ is the class of a simple closed geodesic in
$M$. Then $\Gamma_0=\langle x_0\rangle$ is malnormal in $\Gamma$.
\end{corollary}

In order to apply Proposition \ref{3e} to find the Puk\'anszky invariant of $\vn(\Gamma_0)$, we must determine the size of
$\Gamma_0 \backslash \Gamma /\Gamma_0$. This is done geometrically by considering the diagonal action of $\Gamma$ on the set $$\mathcal F=\{(g_1F_0, g_2F_0) \, :\, g_1, g_2\in \Gamma\}\,.$$
\begin{lemma}\label{F} There is a bijection between $\Gamma_0 \backslash \Gamma /\Gamma_0$
and the set of $\Gamma$-orbits of elements of $\mathcal F$, under the diagonal action.
\end{lemma}

\begin{proof}
The required bijection is the composition of the bijections :
$$
\Gamma(F_0, gF_0) \mapsto \Gamma(\Gamma_0\tilde\xi, g\Gamma_0\tilde\xi)
= \Gamma(\tilde\xi, \Gamma_0g\Gamma_0\tilde\xi)
\mapsto \Gamma_0g\Gamma_0\,,
$$
where $g\in \Gamma$.
The verification of bijectivity is easy, given that $\Gamma$ acts freely on $X$
and that the stabilizer of $F_0$ in $\Gamma$ is $\Gamma_0$.
\end{proof}

\begin{theorem}\label{puk}
Suppose that $M$ has constant negative curvature and that $x_0$ is the class of a simple closed geodesic in
$M$. If $\Gamma_0=\langle x_0\rangle$, then
\begin{itemize}
\item[(a)] $\#(\Gamma_0 \backslash \Gamma /\Gamma_0) = \infty$ ;
\item[(b)] the Puk\'anszky invariant of $\vn(\Gamma_0)$ is $\infty$.
\end{itemize}
\end{theorem}

\begin{proof}
Assuming the curvature is $-1$, the symmetric space $X$ which covers $M$ is a real hyperbolic space of dimension $n\ge 2$, and the maximal flats in
$X$ are geodesics.
By Proposition \ref{3e} and Corollary \ref{mal}, it suffices to prove part (a). By Lemma \ref{F}, this is equivalent to the existence of infinitely many $\Gamma$-orbits of pairs of geodesics in $\mathcal F=\{(g_1F_0, g_2F_0) \, :\, g_1, g_2\in \Gamma\}$,
where the geodesic $F_0$ is the axis of $x_0$. Now the distance in $X$ between the geodesics $gg_1F_0$ and  $gg_2F_0$ is independent of $g\in \Gamma$. It is therefore enough to prove that
there are elements $g\in \Gamma$ for which $d(F_0, gF_0)$ is arbitrarily large.

The unit ball $\{ x\in \bb R^n \, :\, |x|<1\}$, with the appropriate metric will be used as a model for $X$ \cite[1.1]{Ni}, and its boundary sphere $S$ will have its usual metric.
Choose $g\in \Gamma \backslash \Gamma_0$. We show that $d(F_0, g^mF_0)\to\infty$ as $m\to\infty$. The element $g$ is hyperbolic \cite[II.6.3]{BH}. For it cannot be elliptic ($\Gamma$ is torsion free) and it
cannot be parabolic ($\Gamma$ is co-compact).
Therefore $g$ has attracting and repelling fixed points $z^+, z^- \in S$. See, for example \cite[Lemma III.3.3]{Ba}. Now $\{z^+, z^-\}\cap \{z_1, z_2\} =\emptyset$, by Lemma \ref{SST}. It follows that, for $j=1,2$, $g^m z_j \to z^+$ as $m\to \infty$.

\begin{figure}[htbp]
\centerline{
\beginpicture
\setcoordinatesystem units  <0.5in, 0.5in>
\setplotarea x from -3 to 3, y from -2 to 2 
\putrule from  -2 0 to 2 0  
\put {$z_1$}  at  -2.22 0
\put {$z_2$}  at  2.25 0
\put {$_{\bullet}$}  at  0.7 0
\put {$a_m$}  at  0.7 -0.2
\circulararc 360 degrees from 2 0 center at  0 0    
\circulararc 28.5 degrees from 0.78 1.85 center at  8 -1.5    
\put {$z^+$}  at  1  2.1
\put {$z^-$}  at  0  -2.2
\circulararc 163 degrees from .3 1.97 center at  0.78 1.85    
\put {$g^mz_1$}  at  .2  2.25
\put {$g^mz_2$}  at   1.7  1.7
\endpicture
}
\end{figure}

Let $K=\min\{ |z^+-a| \, :\, a\in F_0\} >0$, where $|\cdot |$ is the euclidean norm on $\bb R^n$.
If $0<\epsilon < K$, there exists an integer $m$ such that
$|g^mz_j-z^+|<\epsilon$, $j=1,2$.
There is a unique point $a_m\in F_0$ which is closest (in the hyperbolic metric) to $g^mF_0$.
Let
$$s=d(a_m, g^mF_0)=d(F_0,g^mF_0)\,.$$
 The explicit formula in \cite[Theorem 1.2.1]{Ni} gives

\begin{align*}
\cosh s &=
\frac{2|g^mz_1-a_m||g^mz_2-a_m|}{|g^mz_1-g^mz_2|(1-|a_m|^2)}\\
&\ge
\frac{2(|z^+-a_m|-\epsilon)^2}{2\epsilon}\\
&\ge
\frac{(K-\epsilon)^2}{\epsilon}\,.
\end{align*}
Letting $\epsilon \to 0$ proves the result.
\end{proof}

\end{document}